\documentclass{amsart}[11pt]
\usepackage{myAMSpackages}
\usepackage{Universal}
\usepackage{MyAMSEnvironments}
\usepackage{graphicx,subfigure,epsfig}

\newcommand{\tor}{\text{Tor}}

\newcommand{\tH}{\ensuremath{\text{\tiny $H$\normalsize }}\xspace}
\newcommand{\tV}{\text{\tiny $V$\normalsize }}

\newcommand{\J}{\mathcal{J}}

\newcommand{\Rc}{\ensuremath{\text{Rc}}\xspace}

\renewcommand{\theenumi}{\arabic{enumi}}

\newcommand{\mt}{{\mathfrak{t}}}

\newcommand{\mj}{{\mathfrak{j}}}

\newcommand{\mso}{{\mathfrak{so}}}
\newcommand{\mo}{{\mathfrak{o}}}

\newcommand{\bnabla}{\overline{\nabla}^\lambda}

\newcommand{\tsig}{{\scriptscriptstyle \Sigma}}

\newcommand{\bp}{\begin{pmatrix}}
\newcommand{\ep}{\end{pmatrix}}

\newcommand{\avB}{\Bt{\mathcal{T}}}
\newcommand{\mxB}{\mathcal{T}^+}

\begin{document}
\title{The topology of quaternionic contact manifolds}

\begin{abstract}
We explore the consequences of curvature and torsion on the topology of quaternionic contact manifolds with integrable vertical distribution. We prove a general Myers theorem and establish a Cartan-Hadamard result for almost qc-Einstein manifolds \end{abstract}

\maketitle

\section{Introduction}

In \cite{Biquard}, Biquard introduced quaternionic contact manifolds as the key tool to study the conformal boundaries at infinity of quaternionic K\"abler manifolds. Along with strictly pseudoconvex pseudohermitian manifolds, this class describes a model category of sub-Riemannian manifolds with special holonomy. Subsequently these manifolds themselves have been the objects of extensive study, see for example \cite{AFIV, Duchemin,IMV2,IMV3,IPV, VassilevEinstein,Vassilev1} amongst others. 

Contact quaternionic manifolds  possess a connection adapted to the quaternionic structure \cite{Biquard,Duchemin}. As shown in \cite{VassilevEinstein,Vassilev1} this connection enjoys the remarkable property that its Ricci tensor can decomposed entirely into three torsion components. In this paper, under the assumption that the canonical vertical distribution is integrable, we study the affect of these torsion components on the underlying topology of the manifold. The structure of the paper is as follows. In section 2, we review some basic properties of the Biquard connection and derive some new decompositions of the vertical components of curvature into torsion pieces. In section 3, we show how the Levi-Civita connection for a canonical family of Riemannian metrics can be computed in terms of the Biquard connection and its torsion. In section 4, we establish a general  Bonnet-Myers type theorem.  In section 5, we introduce the category of almost qc-Einstein manifolds, where the horizontal scalar curvature is constant and the horizontal Ricci curvature satisfies
\[ \Rc^\tH (JX,Y) +\Rc^\tH(X,JY)=0\]
for any of the horizontal operators $J$ derived from quaternionic multiplication.  For these manifolds, we employ some techniques of foliation theory to derive Cartan-Hadamard type theorems. In particular, when the horizontal sectional curvatures are non-positive, we show that universal cover is either $\rn{h+3}$ or $\rn{h} \times \sn{3}$ and that the two cases can be distinguished by properties of the torsion.

\section{Basic properties of quaternionic contact manifolds}

In this section, we review the basic definitions associated to quaternionic contact manifolds and establish the basic properties of the Biquard connection. In particular, we derive precise expressions for vertical Ricci and sectional curvatures in terms of torsion components.

\bgD{s2sr}
A step $2$ sub-Riemannian manifold is a smooth manifold $M$ together with a smooth distribution $H \subseteq TM$ of dimension $h$ and a smooth positive definite inner product $g^\tH$ on $H$ such that at every point $[H,H] =TM$.
\enD
If we let $V^* \subseteq T^*M$ consist of the covectors $\xi$ that annihilate $H$, then  we can define a bundle map  $\J \colon V^* \to \text{End}(H)$ by
\[ \aip{ \J(\xi)(X) }{Y}{} = d\xi(X,Y).\]
It is straightforward to verify that $\J(\xi)$ is independent of how $\xi$ is extended to a $1$-form and so is well-defined. The inner product on $H$ allows for pointwise identification of $\text{End}(H)$ with $\mathfrak{gl}_h$, well-defined up to conjugation. If we impose the standard inner product
\[ \aip{A}{B}{} = \text{tr} (B^\top A) \]
on $\mathfrak{gl}_h$, then we can define a canonical inner product on $V^*$ by
\[ \aip{\xi}{\eta}{} = \frac{1}{h} \aip{\J(\xi) }{\J(\eta)}{}.\]
If $V$ is a {\em complement}  to $H$, i.e. a subbundle $V \subseteq TM$ such that $TM = H \oplus V$, then there is a canonical extension of the sub-Riemannian inner product to a Riemannian metric $g$, defined by declaring $H$, $V$ to be orthogonal and using the dual of the inner product of $V^*$ on $V$. For convenience of notation, we denote by $U_\xi $ the dual element in $V$ to $\xi \in V^*$.
\bgD{quaternionic}
A step $2$ sub-Riemannian manifold is quaternionic contact if $\J(V^*)$ is isomorphic to $\mathfrak{sp}_1$ at every point. 
\enD

Thus there is an $\text{SO}(3)$-bundle of triples of unit length forms $\eta^1,\eta^2,\eta^3$ such that for each $a =1,2,3$,
\[ J_{a}^2 =- 1 =J_{123}  \]
where $J_a = \J(\eta^a)$ and if $I=a_1a_2\dots a_k$ then $J_I = J_{a_1} J_{a_2} \dots J_{a_k}$. Furthermore, we can always choose a reduction of the horizontal frame bundle to $\text{Sp}(1)\text{Sp}(h/4)$.  We introduce
\begin{align*}
\mt_0 &= \{ A \in \mso_h \colon [A,\mathfrak{sp}_1]=0\} \cong \mathfrak{sp}_{h/4} \\
\mt &= \mt_0 \oplus \mathfrak{sp}_1
\end{align*}
and denote by $\mt^\perp$ the orthogonal complement in $\mso_h$. It is easily seen that
\[ \mt = \{ A \in \mso_h \colon [A,\mathfrak{sp}_1] \subseteq \mathfrak{sp}_1  \} \cong \mathfrak{sp}_1 \oplus  \mathfrak{sp}_{h/4}.\]

The foundational theorem on quaternionic contact manifolds due to Biquard \cite{Biquard} posits the existence of a connection adapted the the quaternionic structure.

\bgT[Biquard]{Biquard}
If $M$ is quaternionic contact with $h>4$ then there is a unique complement $V$ and connection $\nabla$ such that
\begin{itemize}
\item $H$, $V$, $g$ and $\J$ are parallel
\item $\tor(H,H) \subseteq V$, $\tor(H,V) \subseteq H$
\item For all  $U \in V$, the operator $\tor(U,\cdot) \colon H \to H$ is in $\mt^\perp \oplus \Sigma^h$.
\end{itemize}
\enT
Here $\Sigma_h \subset \mathfrak{gl}_h $ denotes the space of symmetric elements of $\mathfrak{gl}_h$. If $h=4$ then Duchemin \cite{Duchemin} showed that the  same result holds under the additional assumption that there exists a complement where for $a,b \in\{1,2,3\}$, 
\[ d\eta^a(U_b,X) + d\eta^b(U_a,X) =0\]
for all $X \in H$.  The Duchemin condition is equivalent to $1$-normal or $V$-normal in the language of \cite{Hladky4,Hladky5}. A quaternionic contact manifold is often referred to as {\em integrable} if either $h>4$ or it satisfies the Duchemin condition. However, for the sake of brevity, we shall henceforth always assume, unless other stated, that all quaternionic contact manifolds under consideration are integrable.

The property that $\J$ is parallel can easily be shown to be equivalent to the identity
\[  0 \equiv \nabla \tor (X,Y,A). \]
for $X,Y \in H$ and $A \in TM$.

Henceforth, we shall also always assume that $E_1,\dots,E_h$ is an orthonormal frame for $H$ and that $U_1,U_2,U_3$ is an orientable orthonormal frame for $V$ with coframe $\eta^1,\dots,\eta^3$.  The orientability condition is equivalent to $J_{123} = -1$.

For $a \in \{1,2,3\}$, we denote by $a^+,a^-$ the subspaces of $\mathfrak{gl}_h$ that commute and anti-commute with $J_a$ respectively.  It is also easy to see that $\mathfrak{gl}_h = \Psi[3]\oplus\Psi[-1]$ where $\Psi[\lambda]$ is the eigenspace of the invariant Casimir operator $A \mapsto -\sum\limits_a J_a A J_a$. Clearly, the eigenvalue corresponds to the difference in count between the $J_a$'s that commute with $A$ and those that anti-commute. This is also invariant under orthonormal changes of the vertical frame.

 For ease of notation, we denote $\tor(A,B)$ by $T(A,B)$. For $a,b,c \in \{1,2,3\}$ we define torsion operators by  
 \begin{align}
 T_a \colon H \to H, &\qquad T_a X = \tor(U_a,X)\\
 T_a^\tV \colon V \to V & \qquad T_a^\tV =T(U_a,U)_\tV 
 \end{align}
and functions $\tau_{ab}^c$ by
\[  \tau^a_{bc} U = \eta^a \tor(U_b, U_c)\]
 Furthermore, we denote the decomposition of each $T_a$ into symmetric and skew-symmetric components by $T_a = T_a^\tsig + T^{\mo}_a$ and define tensors $T^\tsig = \sum_a \eta^a \otimes T^\tsig_a$, $T^{\mo} = \sum_a \eta^a \otimes T^{\mo}_a$
 
The Biquard connection then has the following useful properties.

 \bgL{easytor2}
For elements $a,b \in \{1,2,3\}$,
\begin{enumerate}
\item $\{J_a,T_b\} \in a^+ \cap b^+$, $\{J_a,T_b\} + \{J_b,T_a\} =0$ and hence $T^\tsig_a \in a^-$.
\item $\text{tr}_\tH \left( T^\tsig_a \right) =0$.
\item $[J_a,T^{\mo}_b] \in  a^- \cap  b^-$, $[J_a,T^{\mo}_b]+[J_b,T^{\mo}_a]=0$ and hence $T^{\mo}_a \in a^+ $.
\item $\tau_{12}^3 = \tau_{23}^1 =\tau_{31}^2 $
\end{enumerate}
where $\{A,B\}$ is the symmetric sum on $\mathfrak{gl}_h$, $\{A,B\} = AB+BA$.
\enL

As a result of the last part, we shall simplify notation by defining the {\em vertical torsion function}
\bgE{tau} \tau =-\tau_{12}^3,\enE
which shall become integral to later results. The reason for the sign choice should also become clear later on.

\pf

The general Bianchi Identity states
 \bgE{Bianchi} 
 \mathscr{C}  \left( R(A,B)C  - \nabla T (A,B,C) +  T(A,T(B,C)) \right)=0 
 \enE
 where $\mathscr{C}$ denotes the cyclic sum. Various projections of this identity, will be of particular use to us: 
if $X,Y,Z$ are sections of $H$, then
 \begin{align} \label{E:B}
  \mathscr{C}   R(X,Y)Z  &= -\mathscr{C} T(X,T(Y,Z)) = \mathscr{C} T(T(X,Y),Z).
 \end{align}
If $X,Y\in H$  then
 \bgE{BV} \mathscr{C}   R(X,Y)U_b    =  -T^\tV_b T(X,Y) +  T(X,T_b Y) - T(Y,T_b X ).
  \enE
If we apply $\eta^a$ to \rfE{BV} we obtain
  \begin{align*}
    \eta^a T_b^\tV T(X,Y) + \eta^a R(X,Y)U_b &=\eta^a(  T(X,T_b Y) - T(Y,T_b X) )  \\
    &= \aip{J_a X }{T_b Y}{} - \aip{J_a Y }{T_b X}{}  \\
    &= \aip{  \left( \{ J_a, T^\tsig_b \}  + [J_a ,T^{\mo}_b ] \right) X}{Y}{}.   
    \end{align*}
  Therefore if we define operators $R_b^a \colon H \to H$ by $\aip{R_b^a X}{Y}{} =\eta^a R(X,Y)U_b $   then 
   \bgE{raw} \{ J_a, T^\tsig_b \}  + [J_a ,T^{\mo}_b ]  = \sum\limits_c \tau^a_{bc} J_c +  R^a_b .\enE
  For any operator $L$, $ \{J_a, L\} \in a^+ $  and $ [J_a,M] \in a^-$ and so 
  \[ \{J_a,T^0_b\} = \tau^a_{ba} J_a + (R^a_b)^{a^+} .\] 
Now if we set $a=b$ and project onto the $a^+$ component we see that $\{J_a,T^\tsig_a\}=0$. Thus $T^\tsig_a \in a^-$ and so for $b\ne a$ we must have  $J_bT^\tsig_a,T^\tsig_aJ_b \in a^+$. From this, we clearly see that $\{J_b,T^\tsig_a\} \in a^+ \cap b^+$.

Next a symmetric application of \rfE{raw} shows that
  \bgE{rawS}  \left( \{ J_a, T^\tsig_b \} +\{ J_b, T^\tsig_a \}  \right) + \left(  [J_a ,T^{\mo}_b ]+ [J_b ,T^{\mo}_a ] \right)  - \sum\limits_c \left(  \tau^a_{bc} +\tau^b_{ac} \right) J_c =0  \enE
  Since $\tau^a_{bb} =0 =\tau^b_{aa}$, the final term is contained in $(a^- + b^-) \cap \mathfrak{sp}_1$. The first term is in $a^+\cap b^+$ as noted earlier. However, as $T^{\mo}_a,T^{\mo}_b \in \mt^\perp$ and
  \[  \aip{ [\mathfrak{sp}_1,\mt^\perp] }{\mathfrak{sp}_1}{} = \aip{ \mt^\perp }{[\mathfrak{sp}_1,\mathfrak{sp}_1]}{} =0.\]
  And so it is easy to see that the middle term of \rfE{rawS} is in  $(a^- + b^-) + \mathfrak{sp}_1^\perp$. Thus the three grouped terms in \rfE{rawS} are mutually orthogonal. Hence each must vanish.
  
  From this it is easy to see that $\tau^a_{bc}$ is fully alternating as a tensor.
  
\epf

The previous lemma can be generalized to contact manifolds based on structures such as the octonions or Clifford algebras that share many properties with the quaternions, but the following result is isolated to the quaternionic case

\bgC{qcB}
For an integrable quaternionic contact manifold,
\[ J_1 T^{\mo}_1 = J_2 T^{\mo}_2=J_3 T^{\mo}_3.\]
Hence for $a,b \in \{1,2,3\}$, $\|T^{\mo}_a\| =\|T^{\mo}_b\|$ and if $a \ne b$, then $\{T^{\mo}_a,T^{\mo}_b\}=0$ and $\aip{T^{\mo}_a X}{T^{\mo}_b X}{}=0$ for all $X \in H$.
\enC

\pf
Since $T^{\mo}_a $ is orthogonal to  $\mt_0$ it must lie within $\Psi[-1]$. As it commutes with $J_a$, it must therefore lie in $b^-$ if $b \ne a$. But  then   \rfL{easytor2} implies that for $a \ne b$, 
\[  J_b T^{\mo}_a = \frac{1}{2} [J_b ,T^{\mo}_a ] = - \frac{1}{2} [J_a ,T^{\mo}_b ]  = -J_a T^{\mo}_b  .\]
The first result follows easily.

We also see that for $a \ne b$, $T^{\mo}_a = -J_{ab} T^{\mo}_b$. Since $J_{ab} \in b^-$ and $T^{\mo}_b \in b^+$, the remaining properties follow also. 

\epf

We define functions $\avB$, $\mxB$  by
\bgE{B} \avB = \| T^{\mo}_1\|, \qquad \mxB = \sup\limits_{|X|=1} | T^{\mo}_1 X| \enE
The content of the previous corollary is that these definitions are independent of the choice of orthonormal frame for $V$, indeed for any unit length $U \in V$, the mean-square and maximum of the eigenvalues of the operator $ iT(U, \cdot) \colon H \to H$ are $\avB^2$ and $\mxB$ respectively.
 
 For convenience, we introduce the invariant trace-free symmetric operators
 \bgE{tfs}
[JT^\tsig] = \sum\limits_a J_aT^\tsig_a, \qquad
[JT^{\mo}] = \sum\limits_a J_a T^{\mo}_a. 
\enE

\bgC{jbn}
The symmetric operator $[JT^{\mo}]$ has norm, $\|[JT^{\mo}] \| =3 \avB$ and largest eigenvalue $3\mxB$.
\enC


For a quaternionic contact manifold it makes sense to split curvature operators into horizontal and vertical pieces. Since the horizontal bundle is fundamental to the definition of a sub-Riemannian manifold, it is desirable to understand how the horizontal components alone affect the geometry and topology of the manifold. Here we shall focus on the Ricci curvatures.

\bgD{Rc}
If $E_1,\dots,E_h$ and $U_1,U_2,U_3$ denote orthonormal bases for $H$ and $V$ respectively, then 
\begin{align*}
\Rc^\tH(A,B) &= \sum\limits_i \aip{R(E_i,A)B}{E_i}{}\\
\Rc^\tV(A,B) &= \sum\limits_a \aip{R(U_a,A)B}{U_a}{}\\
\end{align*}
\enD

The horizontal Ricci curvature on an integrable quaternionic contact  manifold is well-understood. In \cite{Vassilev2} it is shown that that as an operator on $H$, $\Rc^\tH$ is symmetric and has an orthogonal decomposition into $\text{Sp}(1)\text{Sp}(h/4)$ invariant torsion components given by
\bgE{RcH} \Rc^\tH = \left( \frac{h}{4}+2 \right) \tau -\left( \frac{h}{4} +1  \right) [JT^\tsig] - \frac{ h+10}{6}  [JT^{\mo}] \enE
where the difference of a factor of $2$ from \cite{Vassilev2} derives from a minor difference of convention in the relationship between the metric and quaternionic endomorphisms.

  For mixed terms we can use a double cyclic argument (similar to \cite{Hladky5}, Lemma 1) to see
    \begin{align*}
\aip{R(E_i,U)X}{E_i}{} &= \aip{R(E_i,U)X}{E_i}{} - \aip{R(E_i,X)U}{E_i}{} \\
&= \frac{1}{2} \mathscr{C}  \aip{\mathscr{C}R(U,X)E_i}{E_i}{} \\
& = \frac{1}{2} \mathscr{C} \aip{\mathscr{C} T (T(U,X),E_i)}{E_i}{}+\frac{1}{2} \mathscr{C} \aip{\mathscr{C} \nabla T(U,X,E_i)}{E_i}{}  \\
& =\aip{ T(T(X,E_i),U)}{E_i}{} +  \aip{\nabla T(U,X,E_i)}{E_i}{}  \\
& \qquad -\aip{\nabla T (U,E_i,X)}{E_i}{}.
\end{align*}
Since \rfL{easytor2} (2) implies that $\sum\limits_i \aip{T(U,E_i)}{E_i}{} =0$ we can easily see that $\sum\limits_i \aip{\nabla T (U,E_i,X)}{E_i}{}=0$. Hence
    \bgE{Hricmix}
    \Rc^\tH(U,X) = \sum\limits_a \aip{T(U,U_a)}{J_a X}{}  + \sum\limits_ i \aip{\nabla T(U,X,E_i)}{E_i}{}.
    \enE
 Furthermore, if $V$ is integrable the first term of the right vanishes.

We shall primarily be interested in $\Rc^\tV$ as an operator on $V$. It is less well-behaved than $\Rc^\tH$, having both symmetric and skew-symmetric components, however it too can be decomposed into torsion components.

We first note that all purely vertical components of the full curvature tensor of $\nabla$ can actually computed using purely horizontal operators expressible either in terms of curvatures or torsions.

\bgL{RmV} 
For $a,b,c,d\in \{1,2,3\}$,
\[ \aip{R(U_a,U_b)U_c}{U_d}{} = -\frac{2}{h}\aip{ R(U_a,U_b)}{J_{dc}}{}= -\frac{2}{h} \aip{[T_a,T_b]}{J_{dc} }{} \]
where the operator inner products are taken on $\text{End}(H)$.
\enL

\pf
 This begins with an easy computation using the parallel torsion properties. Namely,
  \begin{align*}
\aip{R(U_a,U_b)U_c}{U_d}{} &= \frac{1}{h} \sum\limits_i  \aip{R(U_a,U_b) T(E_i,J_c E_i) }{U_d}{} \\
&= \frac{1}{h}\sum\limits_i   \aip{ T(R(U_a,U_b) E_i,J_c E_i) +T(E_i,R(U_a,U_b) J_c E_i)}{U_d}{}  \\
&= \frac{1}{h} \aip{ J_d R(U_a,U_b)}{J_c}{} + \frac{1}{h} \aip{J_d}{R(U_a,U_b)J_c}{}\\
&= -\frac{2}{h}\aip{ R(U_a,U_b)}{J_{dc}}{}.
\end{align*}
Next  we use the cyclic identities to break the operator $R(U_a,U_b)$ into torsion pieces. For $E \in H$,
\bgE{RVtoH}
\begin{split}
R(U_a,U_b)E &= \left( \mathscr{C} R(U_a,U_b)E \right)_\tH = \mathscr{C} T(T(U_a,U_b),E)_\tH + \mathscr{C} \nabla T (U_a,U_b,E)_\tH \\
&= \tau_{ab}^c T_c E -[T_a,T_b]E + [\nabla_{U_a}, T_b]E  -\aip{\nabla_{U_a} U_b}{U_c}{} T_c E  \\
& \qquad - [\nabla_{U_b}, T_a]E + \aip{\nabla_{U_b} U_a}{U_c}{} T_c E.
\end{split}
\enE
Now 
\begin{align*}
\aip{ [ \nabla_{U_b}, T_a]}{J_c}{} &= U_b \aip{T_a}{J_c}{} - \aip{T_a}{\nabla_{U_b} J_c }{} + \aip{\nabla_{U_b} }{(T^{\mo}_a - T^\tsig_a) J_c}{} \\
&= 0  -\aip{T_a}{ J_{\nabla_{U_b} \eta^c} }{} + \aip{J_c T_a + (T^{\mo}_a - T^\tsig_a) J_c }{\nabla_{U_b}}{} \\
&= \aip{ \{J_c,T^{\mo}_a\} + [J_c,T^\tsig_a] }{\nabla_{U_b}}{}\\
&=0
\end{align*}
as the left hand side of the penultimate line is symmetric but the right skew-symmetric.

Thus, recalling that $J_{dc} \in \mj \oplus \langle I_h \rangle$, and that $T_c$ is always trace free and orthogonal to $\mj$, we see that
 \bgE{RT1} \aip{R(U_a,U_b)U_c}{U_d}{}   = - \frac{2}{h} \aip{[T_a,T_b]}{J_{dc} }{} .\enE
 
\epf

We can pursue this further to compute the vertical sectional curvatures in terms of the torsion operators as follows. 
\bgL{KV}
For $a,b \in \{1,2,3\}$ with $a\ne b$, the vertical sectional curvatures satisfy
\[ \frac{h}{2} K(U_a,U_b)  = \avB^2  -  \big\| (T^\tsig_a)^{b^+} \big\|^2=\avB^2  -  \big\| (T^\tsig_b)^{a^+} \big\|^2.\]
\enL

\pf

We apply \rfE{RT1} with $c=a$ and $b=d$ and note that $J_{ba}$ is either pure trace or skew-symmetric depending on whether $a=b$. Now $T_c$ splits into symmetric and skew-symmetric pieces as $T_c = T^\tsig_c + T^{\mo}_c$.  We then observe that $[T_a,T_b]$ is trace-free and its skew-symmetric component  is $[T^\tsig_a,T^\tsig_b] + [T^{\mo}_a,T^{\mo}_b]$.

Next we see
\begin{align*}
 \aip{ T^\tsig_a T^\tsig_b - T^\tsig_b T^\tsig_a}{J_{ab}}{} &=-2 \aip{J_bT^\tsig_b }{J_a T^\tsig_a}{}
 \end{align*}
 and using the torsion properties we see  that if $a\ne b$
 \begin{align*}
 0 & = \aip{ \{ J_a ,T^\tsig_b\} + \{J_b ,T^\tsig_a\} }{ J_b T^\tsig_a }{} \\
 &= \aip{ -J_b J_a T^\tsig_b -J_bT^\tsig_b J_a}{T^\tsig_a }{} + \aip{ T^\tsig_a -J_b T^\tsig_a J_b}{T^\tsig_a}{}\\
 &= -2 \aip{J_bT^\tsig_b }{J_a T^\tsig_a}{} + \aip{ T^\tsig_a -J_b T^\tsig_a J_b}{T^\tsig_a}{}\\
 &= -2 \aip{J_bT^\tsig_b }{J_a T^\tsig_a}{} + 2 \aip{ (T^\tsig_a)^{b^+} }{T^\tsig_a}{} .
 \end{align*}
 Thus
 \[ \aip{[T^\tsig_a,T^\tsig_b]}{J_{ba}}{} = -2 \big\| (T^\tsig_a)^{b^+}\big\|^2 .\]
 

 
 A similar argument shows that
 \begin{align*}
 \aip{ T^{\mo}_aT^{\mo}_b-T^{\mo}_bT^{\mo}_a }{J_{ba}}{} &=2 \aip{J_b T^{\mo}_b}{J_a T^{\mo}_a}{}=2\avB ^2
\end{align*}
The result follows easily. 
 \epf

The curvature for the Biquard connection does not enjoy all the same symmetries possessed by the Levi-Civita connection. The sectional curvatures do not determine the full curvature tensor via polarization, or  even the Ricci curvatures.. The symmetric portion of $\Rc^\tV$ is of course determined by the sectional curvatures, but unlike $\Rc^\tH$ or the Riemannian Ricci tensor, there is also a skew-symmetric component.

 \bgL{RcV}
 For an  integrable quaternionic contact manifold with integrable vertical complement, if $J_{abc}=-1$ then
 \[ \Rc^\tV(U_a,U_b) = \frac{4}{h} \aip{ T^{\mo}_a }{T^{\mo}_b}{}- \frac{4}{h} \aip{T^\tsig_a }{T^\tsig_b}{} +\frac{1}{2} d\tau(U_c) .\]
 \enL

 \pf
 First, we choose an orthonormal frame $\eta^1,\dots,\eta^3$ for $V^*$ such that $\xi = f\eta^1$. Then
 \[ \Rc^\tV(U_\xi,U_\xi) =f^2 \Rc^\tV(U_1,U_1).\]
 Now since $T^{\mo}_1 \in \Psi[-1]$, it easily follows that
 \[ \Rc^\tV(U_1,U_1) = \sum\limits_a K(U_a,U_1) = \frac{4}{h}\big\| T^{\mo}_1\big\|^2 - \frac{4}{h} \big\| T^\tsig_1\big\|^2.\] 
 The symmetric component of $\Rc^\tV$ can now simply be computed from a polarization argument, so it only remains to find the skew-symmetric part. Now
  \bgE{cyclic}\begin{split}
2 &\aip{R(U_3,U_1)U_2}{U_1}{}  - 2 \aip{R(U_2,U_1)U_3}{U_1}{} \\
&\qquad =\mathscr{C}  \aip{\mathscr{C}R(U_1,U_2)U_3}{U_1}{} \\ 
&\qquad =\mathscr{C} \aip{\mathscr{C} T (T(U_1,U_2),U_3)}{U_1}{} +\mathscr{C} \aip{\mathscr{C} \nabla T(U_1,U_2,U_3)}{U_1}{}  \\
&\qquad= 2 d\tau(U_1) 
\end{split}\enE
Hence
\[  \Rc^\tV(U_2,U_3) -\Rc^\tV(U_3,U_2)  = d\tau(U_1). \]
With an identical argument for the other components, the result follows immediately.
 
 \epf

Somewhat surprisingly, for mixed terms, the vertical Ricci tensor is often better behaved than the horizontal.
 \bgL{MixedRc}
For $U \in V$ and $X \in H$,
 \[ \Rc^\tV(X,U) =-\sum\limits_a \aip{T(U,U_a)}{J_a X}{}  .\]
 If $V$ is integrable then $\Rc^\tV(X,U) =0$.
 \enL
 
 \pf
Here we employ standard results derived from the algebraic Bianchi identity  (see ) for a connection with torsions. Notably
\begin{align*}
\aip{R(U_a,X)U}{U_a}{} &= \aip{R(U_a,X)U}{U_a}{} - \aip{R(U_a,U)X}{U_a}{} \\
&= \frac{1}{2} \mathscr{C}  \aip{\mathscr{C}R(X,U)U_a}{U_a}{} \\
& = \frac{1}{2} \mathscr{C} \aip{\mathscr{C} T (T(X,U),U_a)}{U_a}{}+\frac{1}{2} \mathscr{C} \aip{\mathscr{C} \nabla T(X,U,U_a)}{U_a}{}  \\
& = -\aip{T(U,U_a)}{J_a X}{}  + \aip{\nabla T(U,U_a,X)}{U_a}{} 
\end{align*}

However by skew-symmetry of the purely vertical torsion, we see 
\[ \aip{T(U,U_a)}{U_a}{} =0.\] The result follows easily.
 
 \epf

 \section{Comparison with weighted Levi-Civita connections.}
 
 The extension of the metric on $H$ to a full Riemannian metric $g$  described earlier was canonical only up to a constant scaling factor.   Thus it is natural to consider a family of Riemannian metric defined by
 \[ g^\lambda(A,B) = \aip{A_\tH}{B_\tH}{} + \lambda^2 g(A_\tV,B_\tV) .\]
 While the Biquard connection associated to the quaternionic structure contains more refined information on the geometry of a quaternionic contact manifold, the torsion-free nature of the Levi-Civita connections $\bnabla$ allows for the application of the more deeply developed theory of Riemannian geometry.

 Our first step is to compare the Biquard connection to these weighted Levi-Civita connections. This is technically much simpler if we make the assumption that the vertical distribution is integrable. It should however be noted that this condition appears automatically in many important cases such qc-Einstein manifolds (see \cite{Vassilev2}) and the families involved in the partial solution to the qc-Yamabe problem from the same paper. One consequence of this is that $M$ then admits a foliation $\mathcal{F}^\tV$ whose leaves are $3$-dimensional manifolds everywhere spanned by $V$.
  
  \bgL{ConComp} If the vertical distribution $V$ is integrable then the Levi-Civita derivatives with respect to the metric $g^\lambda$ can be computed for sections $X,Y$ of $H$ as follows:
  \begin{align*}
  \bnabla_X Y &=\nabla_X Y  - \frac{1}{2} \sum\limits_a \aip{J_a X}{Y}{} U_a  - \lambda^{-2} \sum\limits_a  \aip{T^\tsig_a X}{Y}{} U_a,\\
  \bnabla_{U_a} X &= \nabla_{U_a} X + \frac{\lambda^2}{2} J_a X  -T^{\mo}_a X, \\
\bnabla_X U_a  &=  \nabla_X U_a +  \frac{\lambda^2}{2} J_a X +T^\tsig_a X,\\
\bnabla_{U_a}  U_b &= \nabla_a U_b-\frac{1}{2} \sum\limits_c \tau_{ab}^c U_c.
 \end{align*}
 \enL
 Here we get the first indication that the vertical torsion function $\tau$ will play an instrumental role in determining the topology in the leaves of $\mathcal{F}^\tV$.
 
 \pf The proof works by using the standard formulas expressing the Levi-Civita connections in terms of Lie brackets and then decomposes those Lie brackets using the Biquard connection. The computations are  similar for all parts, so we shall prove the first and leave the others to the reader.
 
 We begin by noting that, since the difference between connections is tensorial, it suffices to prove the results for $X,Y$ members of an horizontal orthonormal frame. Then letting $Z$ be another member of the same horizontal frame
 \begin{align*} g^\lambda(\bnabla_X Y,Z) &=-\frac{1}{2} \left( \aip{X}{[Y,Z]}{}  +\aip{Y}{[X,Z]}{} - \aip{Z}{[X,Y]}{} \right)  \\
 &= \aip{\nabla_X Y }{Z}{}
 \end{align*}
 \begin{align*}
 g^\lambda(\bnabla_XY,U_b) &= -\frac{1}{2} \left( \aip{X}{[Y,U_b])}{} +\aip{Y}{[X,U_b]}{} - g^\lambda(U_b,[X,Y]) \right) \\
 &= -\frac{1}{2} \Big( \aip{X}{ -\nabla_{U_b} Y -T_b Y}{} + \aip{Y}{-\nabla_{U_b} X -T_b X}{} \\
 & \qquad  + \lambda^2 g(U_b,T(X,Y) ) \Big)\\
 &= \frac{1}{2} \aip{T^\tsig_b X}{Y}{} -\frac{\lambda^2}{2} \aip{J_b X}{Y}{}.
 \end{align*}
The proof of the first result is completed by then recalling that $\lambda^{-1} U_1 ,\dots, \lambda^{-1} U_3$ is an orthonormal frame for $V$ with respect to $g^\lambda$.

 \epf
 
 
 
 
 With this comparison in hand, we can attend to the laborious task of comparing the curvature tensors.
 
  \bgT{Rm} If $M$ is an integrable quaternionic contact manifold with integrable vertical distribution then the Levi-Civita connections associated to $g^\lambda$ can be computed from the Biquard connection as follows.
  \begin{align*}
 \Bt{\text{Rm}}^\lambda& (X,Y,Z,W) = \text{Rm}(X,Y,Z,W) \\
 & - \lambda^{-2} \sum\limits_a    \big[\frac{\lambda^2}{2} \aip{J_a Y}{Z}{}+  \aip{T^\tsig_a Y}{Z}{}\big] \big[ \frac{\lambda^2}{2} \aip{J_a X}{W}{}+  \aip{T^\tsig_a X}{W}{}\big]  \\
 & + \lambda^{-2} \sum\limits_a    \big[\frac{\lambda^2}{2} \aip{J_a X}{Z}{}+ \aip{T^\tsig_a X}{Z}{}\big] \big[ \frac{\lambda^2}{2} \aip{J_a Y}{W}{}+ \aip{T^\tsig_a Y}{W}{}\big] \\
 & +\frac{\lambda^2}{2}\sum\limits_a \aip{J_a X}{Y}{} \aip{  J_a Z}{W}{} \\
 \Bt{\text{Rm}}^\lambda &(X,Y,U_a,Z) = \aip{ \nabla T^\tsig (U_a,Y,X)}{Z}{}-\aip{ \nabla T^\tsig (U_a,X,Y)}{Z}{}\\
 \Bt{\text{Rm}}^\lambda & (X,U_a,U_a,X) = \aip{\nabla T(X,U_a,U_a)}{X}{} + \frac{\lambda^4}{4} \left|X\right|^2  -\lambda^2 \aip{J_a T^\tsig_a X}{X}{} \\ & \qquad  +\left|T^{\mo}_a X\right|^2 -\left|T_a X\right|^2 \\
  \Bt{\text{Rm}}^\lambda &(U_b,Y,U_a,U_c) = \lambda^2 \text{Rm}(U_b,Y,U_a,U_c)   -\lambda^2 \aip{ \nabla T(U_b,U_a,Y) }{U_c}{}\\
  \Bt{\text{Rm}}^\lambda&(U_a,U_b,U_c,U_d) = \lambda^2 \text{Rm}(U_a,U_b,U_c,U_d) + \frac{\lambda^2}{2}  \left[  d\tau^d_{bc}(U_a) -d \tau^d_{ac}(U_b) \right]\\
 & \qquad + \frac{\lambda^2}{4} \sum\limits_e \left[ \tau^e_{ac}\tau^e_{bd}  - \tau^e_{bc}\tau^e_{ad} -  2 \tau_{ab}^e \tau^e_{cd} \right] 
  \end{align*}
 \enT
 
 \pf
 Again, all parts of the proof are by direct computation. However we can simplify the process considerably by making use of a technical result from \cite{Vassilev2}. Near any point $p \in M$ there always exist orthonormal frames $E_1,\dots,E_h$ and $U_1,\dots,U_3$  such that $(\nabla E_i)_{|p} = 0 =(\nabla U_a)_{|p}$ for all $i,a$.
 
Here, we shall prove the first result and the last and leave the remainder to the reader.  Letting $X,Y,Z,W$ be elements of an orthonormal horizontal frame as above. Then at $p$,
   \begin{align*}
    g^\lambda&( \bnabla_X \bnabla_Y Z, W) =   X \aip{ \bnabla_Y Z}{W}{} - g^\lambda(\bnabla_Y Z, \bnabla_X W) \\
&=  \aip{\nabla_X \nabla_Y Z}{W}{}  \\
& \quad - \lambda^2 \sum\limits_a    \big[\frac{1}{2} \aip{J_a Y}{Z}{}+\lambda^{-2}  \aip{T^\tsig_a Y}{Z}{}\big] \big[ \frac{1}{2} \aip{J_a X}{W}{}+\lambda^{-2}  \aip{T^\tsig_a X}{W}{}\big] 
 \end{align*}
and
 \begin{align*}
 g^\lambda(\bnabla_{[X,Y]}Z,W) &= \sum\limits_a \aip{[X,Y]}{U_a}{} \aip{\bnabla_{U_a} Z}{W}{} \\
 &= -\sum\limits_a \aip{J_a X}{Y}{} \aip{ \frac{\lambda^2}{2} J_a Z}{W}{} 
 \end{align*}
 The first result follows easily.
 
 Now let $U_1,U_2,U_3$ be an orthonormal vertical frame such that each $\nabla U_a $ vanishes at $p$. Then 
  \begin{align*}
 g^\lambda(\bnabla_{U_a} \bnabla_{U_b} U_c,U_d) &=U_a g^\lambda( \bnabla_{U_b} U_c,U_d) -   g^\lambda( \bnabla_{U_b} U_c,\bnabla_{U_a} U_d) \\
 &= \lambda^2 \aip{ \nabla_{U_a} \nabla_{U_b} U_c}{U_d}{} + \frac{\lambda^2}{2}  d\tau_{bc}^d (U_a) - \frac{\lambda^2}{4} \sum\limits_e \tau_{bc}^e \tau_{ad}^e 
 \end{align*}
 \begin{align*}
 g^\lambda(\bnabla_{[U_a,U_b]} U_c,U_d) = - \lambda^2  \sum\limits_e \tau_{ab}^e \aip{\bnabla_{U_e} U_c}{U_d}{} = \frac{\lambda^2}{2} \sum\limits_e \tau_{ab}^e \tau_{ec}^d 
 \end{align*}
 Again the last result follows easily.

 \epf

 From these we get the following simple corollaries.

 \bgC{Sectional} If $X,Y \in H$ are unit length and orthogonal and $a \ne b$ then,
 \begin{align*}
 \Bt{K}^\lambda(X,Y) &= K(X,Y)  -\lambda^{-2}\sum\limits_a  \Big[ \aip{T^\tsig_a X}{X}{} \aip{T^\tsig_a Y}{Y}{}  -\aip{T^\tsig_a X}{Y}{}^2 \Big]  \\
 & \qquad - \frac{3\lambda^2}{4} \sum\limits_a \aip{J_a X}{Y}{}^2\\
 \Bt{K}^\lambda(X,U_a) &= \frac{\lambda^2}{4} - \aip{J_a T^\tsig_a X}{X}{} \\
 & \qquad   + \lambda^{-2} \Big[  \left|T^{\mo}_a X\right|^2 -\left|T_a X\right|^2-\aip{\nabla T^\tsig(U_a, X,U_a)}{X}{}  \Big]\\
 \Bt{K}^\lambda(U_a,U_b) &= \lambda^{-2} \left[ K(U_a,U_b)  + \frac{\tau^2}{4} \right] 
 \end{align*}
 \enC
 \bgC{RcStuff}
 \begin{align*}
 \Bt{\Rc}^\lambda(X,X) &= \Rc^\tH (X,X) -\aip{[JT^\tsig]X}{X}{}  -\frac{3\lambda^2}{2} \left|X\right|^2  \\
 & \quad - \lambda^{-2}\Big[  \aip{ \text{tr}_\tV \nabla T^\tsig(X)}{X}{} +  2  \sum\limits_a \aip{  T^{\mo}_aX}{T^\tsig_aX}{}  \Big]  \\
    \Bt{\Rc}^\lambda(X,U_a) &=\aip{ \text{tr}_\tH \nabla T^\tsig (U_a)}{X}{} \\
   \Bt{\Rc}^\lambda(U_a,U_a) &= \Rc^\tV(U_a,U_a)  +  \frac{\tau^2}{2} + \frac{h \lambda^4}{4} +\|T^{\mo}_a\|^2 -\|T_a\|^2   \\
   &=\frac{\tau^2}{2}+  \frac{h\lambda^4}{4} +  \frac{4\avB^2}{h} - \frac{h+4}{h} \left\| T^\tsig_a \right\|^2  \end{align*} 
 \enC
 
 Since we are assuming that the vertical distribution is integrable, it is natural to ask what affect torsion has on the topology of the leaves of $\mathcal{F}^\tV$. 
 \bgC{RcF}
 The Levi-Civita connection associated to the restriction of any metric $g^\lambda$ to a leaf of the foliation by  $\mathcal{F}^\tV$ has Ricci curvature given by
 \[ \Rc^{L}(U_a,U_a) =\frac{\tau^2}{2} + \Rc^\tV(U_a,U_a)   = \frac{\tau^2}{2} +  \frac{4\avB^2}{h} - \frac{4}{h} \|T_a^\tsig\|^2 \]
 and sectional curvatures
 \[ K^L(U_a,U_b) = \lambda^{-2} \left( \frac{\tau^2}{2} +  \frac{4\avB^2}{h} - \frac{4}{h} \|(T_a^\tsig)^{b+})\|^2 \right).\]
 \enC
 Thus we see that if $T^\tsig$ is small compared to $\tau$ and $\avB$, then the leaves of the foliation will  be compact. Conversely, if $T^\tsig$ is the dominant torsion component, we would have leaves that are covered by $\rn{3}$ and so would either be non-compact or have considerable topology .

 \section{A qc-Bonnet-Myers theorem}
 In this section, we explore conditions that ensure compactness of a quaternionic contact manifold. We shall make the standing assumption that $M$ is complete with respect to one (and hence all) of the metrics $g^\lambda$. It is well known that this is equivalent to $M$ being complete for the underlying sub-Riemannian metric too. The strategy is to choose the scaling factor $\lambda$ carefully and employ the traditional Bonnet-Myers theorem.
 
 Considering \rfC{RcStuff}, the obvious difficulty is that letting $\lambda \to \infty$ produces two opposing effects. The vertical Ricci curvatures become more positive, but the horizontal Ricci curvatures become more negative.
 
 A necessary condition for $M$ to be compact is that the various torsion operators must be bounded. Thus when looking for sufficient conditions, it makes sense to explore the consequences of bounds on particular pieces of torsion.  Thus, throughout the remainder of this section, we shall  suppose that $\ua,\ub,\uc$ are constants satisfying the following pointwise bounds at all points of $M$ and for all  $X \in H$ and $U\in V$:
 \begin{align} -2 \ua  \left| X \right| \, \left|U \right| &\leq   \sum\limits_i \aip{ \text{tr}_\tH \nabla T^\tsig (U)}{X}{}  \\
  -\ub \left|X\right|^2 &\leq  -\aip{ \text{tr}_\tV \nabla T^\tsig (X)}{X}{}-2  \sum\limits_a \aip{T_a^\tsig X }{ T^{\mo}_aX}{} ,\\
\frac{h+4}{h} \left\| T^\tsig(U,\cdot) \right\|^2 & \leq  \left( \uc+ \frac{\tau^2}{2}   +   \frac{4\avB^2}{h}  \right) \left|U\right|^2  .   \end{align}
 It should be remarked that by necessity, we must always have $\ua \geq 0$.  Then
 \begin{align*}
 \Bt{\Rc}^\lambda(X+ \mu U_a,X+ \mu U_a) &\geq  \Rc^\tH (X,X) -\aip{[JT^\tsig]X}{X}{} - \Big[ \frac{\ub}{\lambda^2} +  \frac{3 \lambda^2}{2}+\ua \e  \Big]  \left|X\right|^2  \\
 & \qquad +   \left( \frac{h\lambda^4}{4}   -\uc  - \frac{\ua}{\e}\right)\mu^2  \\
 \end{align*}
 In order to apply the standard Bonnet-Myers theorem, there must be a positive constant $c$, such that
 \[  \Bt{\Rc}^\lambda(X+\mu U_a,X+ \mu U_a) \geq  c \left( |X|^2 +\lambda^2 \mu^2 \right) \]
 regardless of choice of $X$ or $a$.   This clearly requires us to choose
 \[ \e > \frac{ 4\ua}{h\lambda^4 - \uc}.  \]
This then easily leads to the following theorem

\bgT{qcBM} 
Suppose $M$ is a complete, integrable, quaternionic contact manifold with integrable vertical distribution. 
If there exists a constant $\rho_0$ such that
 \bgE{RcC} \Rc^\tH(X,X) - \aip{[JT^\tsig]X}{X}{} \geq \rho_0 \|X\|^2 \enE
with
 \bgE{rhoC}  \rho_0 > \min\limits_{0 < x <x_L}  \left\{   \frac{\ub}{x} +  \frac{3x}{2}+\frac{4\ua^2}{h x^2 -\uc} \right\}, \qquad  x_L = \begin{cases}   \sqrt{\frac{\uc}{h} }, \quad & \uc>0 \\
  \infty, \quad & \uc\leq 0\end{cases} \enE
then $M$ is compact with finite fundamental group.
\enT

We remark that it is mainly the presence of the torsion term $T^\tsig$ that complicates the expressions above. This theorem simplifies considerable if we place additional constraints on this term. Indeed, we can easily establish the following result.

\bgC{TV}
Suppose $M$ is a complete, integrable, quaternionic contact manifold with integrable vertical distribution. 
\begin{itemize}
\item If $\text{tr}_\tH \nabla T^\tsig \equiv 0$ and there is a constant $\rho_0 > \sqrt{\frac{3\ub}{2}}$ such that  \rfE{RcC} holds then $M$ is compact
\item If $T^\tsig\equiv 0$ and there is a constant $\rho_0>0$ such that \rfE{RcC} holds then $M$ is compact.
\end{itemize}
\enC 

It should be remarked however that the combination of integrable vertical distribution and even the first of these conditions is quite restrictive. Indeed it follows from Theorem 4.8 in \cite{Vassilev2} that under such hypothesis $M$ must have constant scalar curvature.



  
  \section{Almost qc-Einstein manifolds}
  
Of particular importance, in the theory of quaternionic contact manifolds is the class of qc-Einstein manifolds for which the Ricci tensor $\Rc^\tH$ is scalar as an operator on $H$.  These manifolds have very special torsion properties due to the following theorem,  shown for $h>4$ in \cite{Vassilev2} and for $h=4$ in \cite{VassilevEinstein}.
 
\bgT{qcEbase}
A integrable quaternionic contact manifold is qc-Einstein if and only if the torsion operators $T^\tsig,B$ vanish identically. Furthermore if $M$ is qc-Einstein then the vertical distribution is integrable and $\tau$ is constant.
\enT

Therefore every qc-Einstein manifold has a vertical foliation $\mathcal{F}^\tV$. Here, motivated in part by the last section, we shall consider a larger class of manifolds that almost meet these conditions.

\bgD{aqcE}
An integrable quaternionic contact manifold is  aqc-Einstein, or almost qc-Einstein, if it has constant scalar curvature ( i.e. $\tau =0$) and the Ricci operator is in $\Psi[3]$ (i.e. commutes with all $J$ operators).

We say that an aqc-Einstein manifold is of noncompact-type if  following three conditions all hold
\begin{itemize}
\item[(NC1)]  $\tau = 0$,
\item[(NC2)]  for every leaf $L$ in $\mathcal{F}^\tV$,  $\inf_L  \avB=0$,
\item[(NC3)]  for every compact leaf $L$ in $\mathcal{F}^\tV$,  $\avB \equiv 0$ on $L$.
\end{itemize}
Otherwise $M$ is of compact-type.
\enD

Theorem  4.8 in \cite{Vassilev2}  implies that $M$ is aqc-Einstein if and only if  $T^\tsig \equiv 0$ and the vertical distribution is integrable. The motivation for the type definitions will come later in \rfL{qcE}.

This means that aqc-Einstein manifolds possess a vertical foliation $\mathcal{F}^\tV$ and a much  simpler form of the Bonnet-Myers theorem.
 
\bgT{qcEp}
 Suppose that $M$ is a complete aqc-Einstein manifold such that 
 \[ \left(\frac{h}{4} +2 \right) \tau - \left( \frac{h+10}{2} \right) \sup\limits_M  \mxB >0.\]  Then $M$ is compact with finite fundamental group.
 
\enT

\pf This an easy consequence of \rfC{TV} together with the decomposition of $\Rc^\tH$ into torsion components..



\epf
 To proceed we shall first need to collect some results from the theory of foliations. The qc-Einstein version of the following result first appeared in \cite{VassilevEinstein}, here we generalize to the aqc-case.
\bgL{qcE}
 Suppose that $M$ is a complete aqc-Einstein quaternionic contact manifold. Then 
  \begin{itemize}
 \item $\mathcal{F}^\tV$ is a Riemannian, bundle-like foliation with respect to each $g^\lambda$ such that each leaf is complete and totally geodesic,
  \item if $M$ is qc-Einstein, then every leaf of $\mathcal{F}^\tV$ is isometric to a complete, $3$-dimensional  space-form with  constant sectional curvature $\tau^2/2\lambda^2$,
 \item if $M$ is of compact-type then every leaf is compact with universal cover $\sn{3}$,
 \item if $M$ is of noncompact-type then every leaf has universal cover $\rn{3}$.
 \end{itemize}
 \enL
 
 \pf
 
 That $\mathcal{F}^\tV$ is totally geodesic and $M$ is fibre-like compatible with $\mathcal{F}^\tV$ follow from the observations that for $X,Y$ sections of $H$ and $U,W$ sections of $V$ respectively,
\[
  \left( \bnabla_U W\right)_\tH=0, \qquad \left(\bnabla_X Y + \nabla_Y X\right)_\tV =0 .\]
  See Lemma 1.2 in \cite{Johnson} and \cite{ONeill} for details. Completeness follows easily from the totally geodesic property.
  
  The second part  follows trivially from the curvature computations   and  \rfL{KV}.   
  
 For the last two parts, we first note Corollary 4 in \cite{Reinhart} stating that all leaves of such a foliation have the same universal covering space.  Now if (NC1) or (NC2) fails, then from \rfC{RcF}  there is a leaf with positive Ricci curvature which must therefore be compact. If (NC3) fails, there is a compact leaf $L$ such that the Ricci curvature on $L$ is quasi-positive. Since $L$ is compact, following the proof of a theorem by Aubin (\cite{aubin}, pp. 398-399) the metric can be deformed on $L$ into one of strictly positive curvature Ricci curvature. In either case, the universal covering space $\tilde{L}$ must then be compact. Now $\tilde{L}$ is compact, simply connected and oriented from the quaternionic identity $J_{123}=-1$. From Poincar\'e duality, we then have $H_1(M)=H_2(M)=0$ and the  Hurewicz theorem then implies that $\pi_2(M) =0$ and $\pi_3(M) =\mathbb{Z}$. By Whitehead's theorem, any generator of $\pi_3(M)$ is then represented by a homotopy equivalence and so  $\tilde{L}$ is homotopy equivalent to $\sn{3}$. Following Perelman's proof of the Poincar\'e conjecture (or Hamilton's partial solution in \cite{Hamilton}), we must have $\tilde{L}$ is diffeomorphic to $\sn{3}$.   As every leaf has the same universal cover, we must have that each leaf is compact.
 
 Now if $M$ is of noncompact-type and $\avB$ vanishes everywhere then $M$ is qc-Einstein with $\tau=0$ and so the universal cover is $\rn{3}$ by the second part.   If $\avB$ is not identically zero, then there must be a noncompact leaf  which admits a point where $\avB>0$. This leaf will then have quasi-positive Ricci curvature by \rfC{RcF}. Its universal cover is then a simply connected, open $3$-manifold with quasi-positive Ricci curvature. A theorem of Zhu, \cite{Zhu}, then implies that the universal cover is $\rn{3}$.

 \epf
 
 While it initially looks difficult to check the type of an aqc-Einstein manifold, it can be reduced to the study of a single leaf.
 
 \bgC{NCT}
 Suppose that $M$ is a complete aqc-manifold with $\tau = 0$ and $p \in M$ satisfies $\avB(p) >0$. Then $M$ is of compact-type if and only if the leaf through $p$ is compact.
 \enC
 
 \pf
Following the  proof of the previous lemma, we note that the leaf through $p$ has quasi-positive Ricci curvature and so has universal cover $\sn{3}$ if it's compact and $\rn{3}$ otherwise. As all leaves have the same universal cover, the result follows immediately.

 \epf
 
We now note the the work of Reinhart on fibre-like foliations, specifically Corollary 3 from \cite{Reinhart}.

\bgT{Rh} If $M$ is a foliated manifold that is complete with respect to bundle-like metric and the foliation is regular,  then $M$ is isometric to a fibre bundle $\pi \colon M \to M^\prime$ where $M^\prime$ is a complete Riemannian manifold and the leaves of $\mathcal{F}$ coincide with the fibres of $\pi$.
 \enT 
 
 It should be remarked here that it is well-known that the regularity condition can be weakened to requiring that leaves of the foliation have trivial leaf holonomy. This is also implied by the simpler condition that each of the leaves is simply connected. 
 
 While this result puts quite strong conditions on the topology of an aqc-Einstein manifold with regular foliation, this regularity condition is often non-trivial to check and indeed often fails. For the remainder of this section, we shall focus on partially replacing this condition with a tensorial, geometric property, namely non-positive horizontal sectional curvatures. As motivation, we note that, if $\mathcal{F}^\tV$ is regular, the horizontal sectional curvatures descend as Riemannian sectional curvatures on the base manifold.
 
 \bgC{qcE}
 If $M$ is a complete aqc-Einstein manifold such that every leaf of $\mathcal{F}^\tV$ has trivial leaf holonomy, then $M$ is a fibre bundle  $\pi\colon M \to M^\prime$ over a Riemannian manifold $M^\prime$ such that for $X,Y \in H$,
 \[ K(X,Y) = K^\prime(\pi_* X,\pi_*Y) \]
 \enC
 
 This is actually independent of which of the metrics $g^\lambda$ we are considering.
 
 \pf The previous theorem establishes that $M$ is a fibre bundle. It follows from standard results on Riemannian submersions that for sections $X,Y$ of $H$ ,
 \[ \pi_* \bnabla_X Y = \nabla^\prime_{\pi_* X}  (\pi_* Y), \qquad \pi_* \bnabla_{U_a} X = \frac{\lambda^2}{2} \pi_* J_a X =  \pi_* \bnabla_X U_a.    \]
 From this it is easy to compute that
 \begin{align*}
  \pi_* \bnabla_X \bnabla_Y Y &=  \nabla^\prime_{\pi_* X}  \nabla^\prime_{\pi_* Y}  (\pi_* Y) \\
   \pi_* \bnabla_Y \bnabla_X Y &=  \nabla^\prime_{\pi_* Y}  \nabla^\prime_{\pi_* X}  (\pi_* Y) -\frac{\lambda^2}{4} \sum\limits_a \aip{J_a X}{Y}{} \pi_*J_a Y \\
    \pi_* \bnabla_{[X,Y]} &= \nabla^\prime_{[\pi_*X,\pi_* Y]} \pi_* Y - \frac{\lambda^2}{2} \sum\limits_a \aip{J_a X}{Y}{}  \pi_* J_a Y
    \end{align*}
    and hence that
    \[ K^\prime (\pi_*X ,\pi_* Y) = \Bt{K}^\lambda(X,Y) + \frac{3\lambda^2}{4} \aip{J_a}{Y}{}^2 = K(X,Y).\]
 \epf
In fact, not only the sectional curvatures descend, but also the quaternionic structure and as noted in \cite{VassilevEinstein} the leaf space $M^\prime$ is actually locally hyper-K\"ahler.  In the case of non-positive sectional curvatures, this can be improved to

 \bgC{aqcE}
 If $M$ is a complete aqc-Einstein manifold such that  $\mathcal{F}^\tV$ is regular the horizontal sectional curvatures are non-positive, then $M$ is diffeomorphic to a fibre-bundle over a manifold $M^\prime$ with universal cover $\rn{h}$.
 \enC

\pf Since $M^\prime$ must have non-positive Riemannian sectional curvatures, has universal cover  $\rn{h}$  by the standard Cartan-Hadamard theorem.

\epf

In fact, we can improve this result by using Hebda's generalization of the Cartan-Hadamard theorem (\cite{Hebda} Theorem 2). This states that for a Riemannian foliation on $M$ with complete bundle-like metric, if the transverse sectional curvatures are non-positive then the universal cover of $M$ fibers over a complete simply-connected Riemannian manifold $M^\prime$ with non-positive sectional curvature and the fibres are the universal cover of the leaves (which we recall must all be the same).  Thus we have established the following aqc-Cartan-Hadamard theorem.

\bgT{qcEm}
Suppose that $M$ is a complete aqc-Einstein manifold such that all the horizontal sectional curvatures  for the Biquard connection are non-positive. Then the universal cover of $M$ is diffeomorphic to $\rn{h} \times \sn{3}$ if $M$ is of compact-type and diffeomorphic to $\rn{h+3}$ if $M$ is of non-compact type.

\enT

It should be remarked here that $[JT^{\mo}]$ is symmetric and trace-free so must have non-negative eigenvalues. The non-positive sectional curvatures imply non-positive Ricci curvature and so require $\tau \leq 0$. It should also be mentioned that the remarkable fact that the torsion terms determine the Ricci tensor does not appear to extend to the sectional curvatures and so the non-positive sectional curvature condition is unlikely to be redundant. 

We conclude with the observation that while it would be useful to remove the restrictive aqc-Einstein condition here, the vertical spaces in the presence of torsion appear to be much wilder. Even in the case where there is a vertical foliation, the foliation would not be Riemannian and we would not expect any analogue of the fibration results. In particular, it would seem unlikely there would be a similarly simple classification of the universal covering spaces under the non-positive sectional curvature condition.

\bibliographystyle{plain}
\bibliography{References}

\begin{thebibliography}{10}

\bibitem{aubin}
T.~Aubin.
\newblock {M\'etriques riemanniennes et courbure}.
\newblock {\em {J. Diff. Geom.}}, 4:383--424, 1970.

\bibitem{Biquard}
O~Biquard.
\newblock { M\'etriques d'Einstein asymptotiquement sym\'etriques}.
\newblock {\em {Ast\'erisque}}, 265, 2000.

\bibitem{AFIV}
L.~de~Andres, M.~Fernandez, S.~Ivanov, L.~Ugarte, and D.~Vassilev.
\newblock {Quaternionic Kaehler and Spin(7) metrics arising from quaternionic
  contact Einstein structures}.
\newblock {\em {Annali di matematica Pura ed Applicata}}, 2012.
\newblock DOI 10.1007/s10231-012-0276-8.

\bibitem{Duchemin}
D~Duchemin.
\newblock { Quaternionic contact structures in dimension 7}.
\newblock {\em {Ann. Inst. Fourier}}, 56(3):851--885, 2006.

\bibitem{Hamilton}
R.S. Hamilton.
\newblock {Three-manifolds with positive Ricci curvature}.
\newblock {\em {J. Diff. Geom}}, 17(2):255--306, 1982.

\bibitem{Hebda}
J.J. Hebda.
\newblock {Curvature and focal points in Riemannian foliations}.
\newblock {\em {Ind. U. Math. J.}}, 35(2):321--331, 1986.

\bibitem{Hladky4}
R.K. Hladky.
\newblock {Connections and curvature in sub-Riemannian geometry}.
\newblock {\em Houston J. Math.}, 38(4):1107--1134, 2012.

\bibitem{Hladky5}
R.K. Hladky.
\newblock {Bounds for the first eigenvalue of the horizontal Laplacian on
  postively curved sub-Riemannian manifolds}.
\newblock {\em Geom. Dedicata.}, 164:155--177, 2013.

\bibitem{IMV2}
S.~Ivanov, I.~Minchev, and D.~Vassilev.
\newblock {Extremals for the Sobolev inequality on the seven dimensional
  quaternionic Heisenberg group and the quaternionic contact Yamabe problem}.
\newblock {\em {J. Eur. Math. Soc.}}, 12:1041--1067, 2010.

\bibitem{IMV3}
S.~Ivanov, I.~Minchev, and D.~Vassilev.
\newblock {The optimal constant in the L2 Folland-Stein inequality on the
  quaternionic Heisenberg group}.
\newblock {\em {Ann. Sc. Norm. Super Pisa Cl. Sci.}}, XI(5):635--652, 2012.

\bibitem{IPV}
S.~Ivanov, A.~Petkov, and D.~Vassilev.
\newblock {The sharp lower bound of the first eigenvalue of the sub-Laplacian
  on a quaternionic contact manifold}.
\newblock {\em {J. Geom. Anal.}}
\newblock to appear, arXiv:1112.0779.

\bibitem{Johnson}
D.L. Johnson and L.B. Whitt.
\newblock {Totally geodesic foliations}.
\newblock {\em {J. Diff. Geo.}}, 15(2):225--235, 1980.

\bibitem{ONeill}
B.~O'Neill.
\newblock {The fundamental equations of a submersion}.
\newblock {\em {Mich. Math. J.}}, 13:459--460, 1966.

\bibitem{Reinhart}
B.L. Reinhart.
\newblock {Foliated manifolds with bundle-like metrics}.
\newblock {\em {Ann. Math. (2)}}, 69(1):119--132, 1959.

\bibitem{VassilevEinstein}
D.~Vassilev and S.~Ivanov.
\newblock { Quaternionic contact Einstein manifolds}.
\newblock {\em {}}.
\newblock arXiv:1306.0474v1.

\bibitem{Vassilev1}
D~Vassilev and S~Ivanov.
\newblock { The sharp lower bound of the first eigenvalue of the sub-Laplacian
  on a quaternionic contact manifold}.
\newblock {\em {J. Geo. Anal.}}
\newblock to appear, arXiv:1112.0779.

\bibitem{Vassilev2}
D~Vassilev, S~Ivanov, and I~Minchev.
\newblock {Quaternionic contact Einstein structures and the quaternionic
  contact Yamabe problem,}.
\newblock {\em {Mem. AMS}}.
\newblock to appear.

\bibitem{Zhu}
S.-H. Zhu.
\newblock {On open three-manifolds of quasi-positive Ricci curvature}.
\newblock {\em {Proc. Ams. Math. Soc.}}, 120(2):569--572, 1994.

\end{thebibliography}
 \end{document}